\documentclass[review]{elsarticle}

\usepackage{subfigure}
\usepackage{float}
\usepackage{amsfonts}
\usepackage{verbatim}
\usepackage{amssymb}
\usepackage{amsmath}

\usepackage{algorithm}
\usepackage{algorithmic}

\usepackage{lineno,hyperref}
\newtheorem{lemma}{Lemma}

\modulolinenumbers[5]

\journal{Indagationes mathematicae-new series}

\begin{document}

\begin{frontmatter}

\title{The ternary Goldbach-Vinogradov theorem with almost equal Piatetski-Shapiro primes}

\author[mymainaddress]{Yanbo Song\corref{mycorrespondingauthor}}
\cortext[mycorrespondingauthor]{Corresponding author}
\ead{syb888@stu.xjtu.edu.cn}

\address[mymainaddress]{School of Mathematics and Statistics, Xi'an Jiaotong University, Xi'an, 710049, China}

\begin{abstract}
In this paper, we proved a theorem that every large enough odd number can be represented as the sum of three almost equal Piatetski-Shapiro primes.
\end{abstract}


\end{frontmatter}

\linenumbers

\section{Introduction}
In 1937, it was proved by I.M. Vinogradov that every large enough odd number could be represented as the sum of three primes. Today we usually called this theorem the ternary Goldbach-Vinogradov theorem. Later many scholars studied the problem that large enough odd number can be represented as the sum of three primes. In \cite{9,10}, C. D. Pan and C. B. Pan showed that every large enough odd number $n$ the Diophantine equation
\begin{align}
n=p_1+p_2+p_3,\ \ \ |p_i-\frac{n}{3}|\leq y,\ \ \ 1\leq i\leq3,
\end{align}
with prime variable is solvable for $y=n^{2/3}(\log n)^c$, where $c$ is a positive constant.

The set of Piatetski-Shapiro primes is of the form
\begin{eqnarray*}
\mathcal{A}_{\gamma}=\{prime\ p|p=[n^{1/\gamma}],\ n\geq1\}.
\end{eqnarray*}
For more information about Piatetski-Shapiro primes, we refer to \cite{4,5,8,12,13,14}.

In 1992, A. Balog and J. Friedlander \cite{1} studied the ternary Goldbach-Vinogradov theorem restricted to the set of Piatetski-Shapiro primes. They proved the following result:
\begin{eqnarray}
\mathop{\mathop{\sum\limits_{p_1+p_2+p_3=n}}}_{p_i\in\mathcal{A}_{\gamma}}\frac{1}{\gamma} p_{1}^{1-\gamma}\frac{1}{\gamma}  p_{2}^{1-\gamma}\frac{1}{\gamma} p_{3}^{1-\gamma}
=\mathfrak{S}_3(n)\frac{3n^2}{\log^3 n}+O\Bigg(\frac{n^2}{\log^4 n}\Bigg),
\end{eqnarray}
where $\mathfrak{S}_3(n)=\prod\limits_{p|n}(1-\frac{1}{(p-1)^2})\prod\limits_{p\nmid n}(1+\frac{1}{(p-1)^3})$ and $20/21<\gamma<1$.
Later C.H. Jia \cite{6} extended the range of $\gamma$ to $15/16<\gamma<1$,
and A.Kumchev \cite{7} improved the result to $50/53<\gamma<1$.

In this paper, we combined these two ideas to study whether (1.1) can hold restricted to the set of Piatetski-Shapiro primes.
Specifically, we proved the following theorem that gives an affirmative answer.\\
\textbf{Theorem} 1. \
For every large enough odd number $n$ the following equation
\begin{eqnarray}
n=p_1+p_2+p_3,\ \ \ |p_i-\frac{n}{3}|\leq y,\ \ \ p_i\in\mathcal{A}_{\gamma}\ \ \ 1\leq i\leq3,
\end{eqnarray}
with prime variable is solvable for $y>n^{(52-45\gamma)/8}$ and $44/45<\gamma<1$.
\section{Notations}

Throughout the paper, the letters $m$ and $n$ are natural numbers, and $p$ always denotes a prime number. Furthermore, $c>1$ is a fixed real number and we put $\gamma=\frac{1}{c}$. And we denote the $\mathcal{A}_c$ the set of the Piatetski-Shapiro primes.

We shall frequently use $\varepsilon$ with a slight abuse of notation to mean a small positive number,
possibly different from line to line.

Given a real number $x$, we write $e(x)=e^{2\pi ix}$, $\{x\}$ for the fractional part of $x$, $[x]$ for the greatest integer not exceeding $x$. We write $\mathcal{L}=\log N$.

We recall that for functions $F$ and real nonnegative $G$ the notations $F\ll G$ and $F=O(G)$ are equivalent to the statement that the inequality $|F|\leq\alpha G$ holds for some constant $\alpha>0$. We also write $F\asymp G$ to indicate that $F\ll G$ and $G\ll F$.

For any function $f$, we put
\begin{eqnarray*}
\Delta f(x)=f(-(x+1)^{\gamma})-f(-x^{\gamma}),\ \ \ (x>0).
\end{eqnarray*}
\section{Some lemmas}
\begin{lemma}
Put $\psi(x)=x-[x]-1/2$. Then, there are numbers $a_h$, $b_h$ such that\\
\begin{eqnarray*}
\left|\psi(x)-\sum\limits_{0<|h|\leq H}a_he(xh)\right|\leq\sum\limits_{|h|\leq H}b_he(xh),
\end{eqnarray*}
where $a_h\ll \frac{1}{|h|}$, $b_h\ll\frac{1}{H}$. Moreover, we let $\psi'(x)=\sum\limits_{0<|h|\leq H}a_he(xh)$, and $\psi''(x)=\psi(x)-\psi'(x)$.\\
\textbf{Proof}.\ For the proof of this lemma, you can see \cite{2}, Appendix.
\end{lemma}
\begin{lemma}
Let $u$, $v\geq1$ be real numbers. For any $n>u$, we have
\begin{eqnarray*}
\Lambda(n)=\mathop{\mathop{\sum_{cd=n}}}_{d\leq v}\log c\mu(d)-\sum\limits_{k|n}a(k)-\mathop{\mathop{\mathop{\sum_{kc=n}}}_{k>1}}_{c>u}\Lambda(c)b(k),
\end{eqnarray*}
where
\begin{eqnarray*}
a(k)=\mathop{\mathop{\mathop{\sum_{cd=k}}_{c\leq u}}}_{d\leq v}\Lambda(c)\mu(d),\ \ b(k)=\mathop{\mathop{\sum_{d|k}}}_{d\leq v}\mu(d).\\
\end{eqnarray*}
\textbf{Proof}.\ This is Prop. 13.4 in \cite{4}.
\end{lemma}
\begin{lemma}
Let
\begin{eqnarray*}
g(\alpha)=\sum\limits_{|p-\frac{n}{3}|<y}e(\alpha p),
\end{eqnarray*}
we have
\begin{eqnarray*}
\int_0^1g(\alpha)^3e(-n\alpha)d\alpha=\mathfrak{S}_3(n)\frac{3y^2}{\log^3 n}+O(\frac{y^2}{\log^4 n}),
\end{eqnarray*}
where $\mathfrak{S}_3(n)=\prod\limits_{p|n}(1-\frac{1}{(p-1)^2})\prod\limits_{p\nmid n}(1+\frac{1}{(p-1)^3})$.\\
\textbf{Proof}.\ This is Theorem 3 in \cite{11}.
\end{lemma}

\begin{lemma}
We put
\begin{eqnarray*}
f(\alpha)=\mathop{\mathop{\sum\limits_{|p-\frac{n}{3}|<y}}}_{p\in \mathcal{A}_1/\gamma}\frac{1}{\gamma}p^{1-\gamma}e(\alpha p).
\end{eqnarray*}
And $g(\alpha)$ are defined in the previous lemma. For any $\alpha$, there are coprime integers $a$ and $q$ with $0\leq a\leq q\leq\log n$ such that $|q\alpha-a|<\log^{-1}n$. and $y>n^{(52-45\gamma)/8}$, $44/45<\gamma<1$, Then we have
\begin{eqnarray*}
f(\alpha)=g(\alpha)+O(n^{21/31-14\gamma/31+\varepsilon}y^{23/31}).
\end{eqnarray*}
\textbf{Proof.}
Let $\chi(n)$ be the characteristic function  of the set $\mathcal{A}_c$, then we have
\begin{eqnarray*}
\chi(n)=[-n^{\gamma}]-[-(n+1)^{\gamma}]=\gamma n^{\gamma-1}+(\psi(-(n+1)^{\gamma})-\psi(-n^{\gamma}))+O(n^{\gamma-2}).
\end{eqnarray*}
Thus,
\begin{eqnarray*}
f(\alpha)=g(\alpha)+\sum\limits_{|p-\frac{n}{3}|<y}cp^{1-\gamma}e(\alpha p)(\psi(-(p+1)^{\gamma})-\psi(-p^{\gamma}))+O(1/\log x).
\end{eqnarray*}
It is enough to estimate,
\begin{eqnarray*}
&&\sum\limits_{|p-\frac{n}{3}|<y}cp^{1-\gamma}e(\alpha p)(\psi(-(p+1)^{\gamma})-\psi(-p^{\gamma}))\\
&=&\sum\limits_{|p-\frac{n}{3}|<y}cp^{1-\gamma}e(\alpha p)\Delta\psi'(p)+O\left(\sum\limits_{|p-\frac{n}{3}|<y}cp^{1-\gamma}e(\alpha p)\Delta\psi''(p)\right)\\
&:=&S_1+S_2.
\end{eqnarray*}
By Lemma 3.1, we have
\begin{eqnarray*}
S_2\ll H^{-1}yn^{1-\gamma}+H^{-1}\sum\limits_{1\leq h<H}\left|\sum\limits_{|x-n/3|<y}x^{1-\gamma}e(hx^{\gamma})\right|.
\end{eqnarray*}
Partial summation yields
\begin{eqnarray*}
\sum\limits_{|x-n/3|<y}x^{1-\gamma}e(hx^{\gamma})\ll n^{1-\gamma}\max\limits_{t<y}\left|\sum\limits_{|x-n/3|<t}e(hx^{\gamma})\right|.
\end{eqnarray*}
Using the exponent pair $(1/2,1/2)$ (\cite{2}, eqn. 3.3.4), we have
\begin{eqnarray*}
\sum\limits_{|x-n/3|<t}e(hx^{\gamma})\ll h^{1/2}n^{\gamma/2}+h^{-1}n^{1-\gamma}.
\end{eqnarray*}
Concluding the above estimate, we have
\begin{eqnarray*}
S_2\ll H^{1/2}n^{1-\gamma/2}+H^{-1}yn^{1-\gamma}.
\end{eqnarray*}
For $S_1$, we first observe that
\begin{eqnarray*}
S_1\ll \frac{1}{\log n}\max\limits_{t<y}\left|\sum\limits_{|x-n/3|<t}cx^{1-\gamma}e(\alpha x)\Lambda(n)\Delta\psi'(x)\right|+O(n^{3/2-\gamma}).
\end{eqnarray*}
Noting that $a_h\ll |h|^{-1}$, we have
\begin{eqnarray*}
\Delta\psi'(x)=\sum\limits_{1\leq|h|\leq H}a_he(hx^{\gamma})(e(h(x+1)^{\gamma}-hx^{\gamma})-1).
\end{eqnarray*}
So we by partial summation, we have
\begin{eqnarray*}
\sum\limits_{|x-n/3|<t}cx^{1-\gamma}e(\alpha x)\Lambda(n)\Delta\psi'(x)\ll \max\limits_{z<t}\sum\limits_{1\leq|h|\leq H}\left|\sum\limits_{|x-n/3|<z}\Lambda(x)e(\alpha n+hx^{\gamma})\right|.
\end{eqnarray*}
We have shown so far
\begin{eqnarray*}
S_1\ll \frac{1}{\log n}\max\limits_{z<y}\sum\limits_{1\leq|h|\leq H}\left|\sum\limits_{|x-n/3|<z}\Lambda(x)e(\alpha n+hx^{\gamma})\right|+O(n^{3/2-\gamma}).
\end{eqnarray*}
Because for any $\alpha$, there are coprime integers $a$ and $q$ with $0\leq a\leq q\leq\log n$ such that $|q\alpha-a|<\log^{-1}n$. Then we have
\begin{eqnarray*}
\sum\limits_{|x-n/3|<z}\Lambda(x)e(\alpha n+hx^{\gamma})=q^{-1}\sum\limits_{-q/2<b\leq q/2}S(a+b,q)\sum\limits_{|x-n/3|<z}\Lambda(x)e(R_b(x)),
\end{eqnarray*}
where $R_b(t)=\beta t+ht^{\gamma}-bt/q$, $\beta=\alpha-a/q$ and
\begin{eqnarray*}
S(a,q)=\sum\limits_{m=1}^{q}e(\frac{am}{q}).
\end{eqnarray*}
Apply Lemma 3.2 with $uv\leq n$ we have
\begin{eqnarray*}
\sum\limits_{|x-n/3|<z}\Lambda(x)e(R_b(x))=\Sigma_1+\Sigma_2+\Sigma_3,
\end{eqnarray*}
where
\begin{align*}
\Sigma_1&=\sum\limits_{r\leq u}\mu(r)\sum\limits_{|rs-n/3|<z}e(R_b(rs))\log s-\mathop{\mathop{\sum\limits_{|rs-n/3|<z}}}_{r\leq u}(\mathop{\mathop{\mathop{\sum\limits_{r=wt}}}_{w\leq u}}_{t\leq v}\mu(w)\Lambda(t))e(R_b(rs)),\\
\Sigma_2&=-\mathop{\mathop{\mathop{\sum\limits_{|rs-n/3|<z}}}_{s>v}}_{r>u}\Lambda(s)\mathop{\mathop{\sum\limits_{d|r}}}_{d\leq u}\mu(d)e(R_b(rs)),\\
\Sigma_3&=-\mathop{\mathop{\sum\limits_{|rs-n/3|<z}}}_{u<r\leq uv}(\mathop{\mathop{\mathop{\sum\limits_{r=wt}}}_{w\leq u}}_{t\leq v}\mu(w)\Lambda(t))e(R_b(rs)).
\end{align*}
By partial summation
\begin{eqnarray*}
\Sigma_1\ll \log x\sum\limits_{r\leq u}\max\limits_{z\leq y}\left|\sum\limits_{|rs-n/3|<z}e(R_b(rs))\right|.
\end{eqnarray*}
Let $f(s)=R_b(rs)$, and noting that $rs\sim n$ and $f''(s)\asymp r^2n^{\gamma-2}|h|$. And applying the van der Corput's estimate in \cite{3} (Cor. 8.13) to $f(s)$ in the inner sum, we have
\begin{eqnarray*}
\Sigma_1\ll \log^2x(u|h|^{1/2}x^{\gamma/2-1}y+x^{1-\gamma/2}|h|^{-1/2}).
\end{eqnarray*}
For $\Sigma_2$, we use dyadic division we can write
\begin{eqnarray*}
\Sigma_2\ll x^{\varepsilon}\sum\limits_{R,S}|T(R,S)|,
\end{eqnarray*}
where
\begin{eqnarray*}
T(R,S)=\sum\limits_{r\sim R}a_r\mathop{\mathop{\sum\limits_{s\sim S}}}_{|rs-n/3|<z}b_se(R_b(rs))
\end{eqnarray*}
with $R>u$, $s>v$, $RS\asymp x$, and $|b_s|$,$|a_r|\leq1$. By Lemma 2.5 in \cite{2} , we have
\begin{eqnarray*}
T(R,S)^2\ll \frac{(RS)^2}{L}+\frac{RS^2}{L}\sum\limits_{1\leq|l|\leq L}\max\limits_{S<s,s+l\leq2S}\sum\limits_{r\in I}e(R_b(r(s+l))-R_b(rs)),
\end{eqnarray*}
where $1\leq L\leq S$ is to be chosen later, and $I\subseteq(R,2R]$ is an interval determined by the conditions $r\sim R$, $|rs-n/3|<z$, $|r(s+l)-n/3|<z$. Let $g(r)=R_b(r(s+l))-R_b(rs)$, so we have $g''(r)\asymp x^{\gamma-1}|hl|R^{-1}$. By Cor. 8.13 in \cite{3}, we have
\begin{eqnarray*}
\sum\limits_{r\in I}e(R_b(r(s+l))-R_b(rs))\ll R^{1/2}x^{\gamma/2-3/2}y|hl|^{1/2}+x^{1/2-\gamma/2}R^{1/2}|hl|^{-1/2}.
\end{eqnarray*}
So we have
\begin{align*}
T(R,S)\ll&xL^{-1/2}+R^{3/4}SL^{1/4}x^{\gamma/4-3/4}y^{1/2}|h|^{1/4}\\
&+ R^{3/4}SL^{-1/4}x^{1/4-\gamma/4}|h|^{-1/4}.
\end{align*}
And using Lemma 2.4 in \cite{2} to choose $L$ optimally we obtain
\begin{align*}
T(R,S)\ll& R^{-1/6}n^{(\gamma+3)/6}y^{1/3}|h|^{1/6}+n^{1-\gamma/4}|h|^{-1/4}\\
&+R^{1/2}n^{1/2}+R^{-1/4}n^{(\gamma+1)/4}y^{1/2}|h|^{1/4}+n^{3/4}y^{1/4}R^{-1/4}.
\end{align*}
So we have
\begin{align*}
x^{-\varepsilon}\Sigma_2\ll&nv^{-1/2}+n^{1-\gamma/4}|h|^{-1/4}+n^{3/4}y^{1/4}u^{-1/4}\\
&+u^{-1/6}n^{(\gamma+3)/6}y^{1/3}|h|^{1/6}+u^{-1/4}n^{(\gamma+1)/4}y^{1/2}|h|^{1/4}.
\end{align*}
Similarly we have
\begin{align*}
x^{-\varepsilon}\Sigma_3\ll&n^{1/2}(uv)^{1/2}++n^{1-\gamma/4}|h|^{-1/4}+n^{3/4}y^{1/4}u^{-1/4}\\
&+u^{-1/6}n^{(\gamma+3)/6}y^{1/3}|h|^{1/6}+u^{-1/4}n^{(\gamma+1)/4}y^{1/2}|h|^{1/4}.
\end{align*}
So we have
\begin{align*}
x^{-\varepsilon}(\Sigma_2+\Sigma_3)\ll&nv^{-1/2}+n^{1/2}(uv)^{1/2}+n^{1-\gamma/4}|h|^{-1/4}+n^{3/4}y^{1/4}u^{-1/4}\\
&+u^{-1/6}n^{(\gamma+3)/6}y^{1/3}|h|^{1/6}+u^{-1/4}n^{(\gamma+1)/4}y^{1/2}|h|^{1/4}.
\end{align*}
Choosing $v=(n/u)^{1/2}$, we have
\begin{align*}
x^{-\varepsilon}(\Sigma_1+\Sigma_2+\Sigma_3)\ll&u^{1/4}x^{3/4}+n^{1-\gamma/4}|h|^{-1/4}+n^{3/4}y^{1/4}u^{-1/4}\\
&+u^{-1/6}n^{(\gamma+3)/6}y^{1/3}|h|^{1/6}+u^{-1/4}n^{(\gamma+1)/4}y^{1/2}|h|^{1/4}\\
&+u|h|^{1/2}x^{\gamma/2-1}y.
\end{align*}
And using Lemma 2.4 in \cite{2} to choose $u$ optimally we obtain
\begin{align*}
x^{-\varepsilon}(\Sigma_1+\Sigma_2+\Sigma_3)\ll&n^{1-\gamma/4}|h|^{-1/4}+x^{3/4}y^{1/8}+|h|^{1/2}x^{\gamma/2-1}y+n^{(\gamma+2)/6}y^{1/3}|h|^{1/6}\\
&+n^{\gamma/4}y^{1/2}|h|^{1/4}+n^{\gamma/8+1/2}y^{1/4}|h|^{1/8}+n^{(\gamma+6)/10}y^{1/5}|h|^{1/10}\\&+n^{\gamma/10+2/5}y^{2/5}|h|^{1/10}
+n^{3\gamma/10}y^{3/5}|h|^{3/10}+n^{3\gamma/14+2/7}y^{3/7}|h|^{3/14}.
\end{align*}
Noting the identity
\begin{eqnarray*}
\sum_{a=0}^{p-1}e\left(\frac{ma}{p}\right)=\left\{
\begin{array}{ll}\displaystyle p, &\textrm{ if $p\mid m$;} \\ \displaystyle
0, &\textrm{ if $p\nmid m$ }\end{array}\right.
\end{eqnarray*}
and sum over $h$, we have
\begin{align*}
S_1\ll& n^{3/2-\gamma}+n^{\varepsilon}(n^{1-\gamma/4}H^{3/4}+n^{(\gamma+4)/8}y^{1/4}H^{9/8}+n^{(3\gamma+4)/14}y^{3/7}H^{17/14}\\&+n^{3/4}y^{1/8}H
+n^{\gamma/2-1}yH^{3/2}+n^{3\gamma/10}y^{3/5}H^{13/10}\\&+n^{(\gamma+6)/10}y^{1/5}H^{11/10}+n^{(\gamma+2)/6}y^{1/3}H^{7/6}+n^{\gamma}y^{1/2}H^{5/4}.
\end{align*}
Combining with $S_2$, and using Lemma 2.4 in \cite{2} to choose $H$ optimally, we have
\begin{align*}
f(\alpha)-g(\alpha)\ll&n^{3/2-\gamma}+n^{\varepsilon}(n^{1-\gamma/4}+n^{(\gamma+2)/6}y^{1/3}+n^{\gamma/4}y^{1/2}+n^{3/4}y^{1/8}+n^{\gamma/8+1/2}y^{1/4}\\
&+n^{\gamma/10+3/5}y^{1/5}+n^{3\gamma/10}y^{3/5}+n^{3\gamma/14+2/7}y^{3/7}+n^{1-\gamma/3}y^{1/3}\\&+n^{1-4\gamma/7}y^{3/7}+n^{7/8-\gamma/2}y^{1/2}
+n^{1/5-2\gamma/5}y+n^{9/13-6\gamma/13}y^{9/13}\\&+n^{5/9-4\gamma/9}y^{7/9}+n^{7/8-\gamma/2}y^{9/16}+n^{13/17-8\gamma/17}y^{11/17}\\
&+n^{17/21-10\gamma/21}y^{13/21}+n^{13/23-10\gamma/23}y^{19/23}+n^{21/31-14\gamma/31}y^{23/31}).
\end{align*}
The proof of the lemma is complete by noting that $n^{21/31-14\gamma/31+\varepsilon}y^{23/31}$ dominates the other terms.
\end{lemma}
\begin{lemma}
Let $f(\alpha)$, $g(\alpha)$ are defined in the previous lemma, then we have
\begin{eqnarray}
&&\int_0^1|f(\alpha)|^2d\alpha\ll n^{1-\gamma}y,\\
&&\int_0^1|g(\alpha)|^2d\alpha\ll y.
\end{eqnarray}
\textbf{Proof.}
We only prove (3.1), since the proof of (3.2) is similar but simpler.
\begin{eqnarray*}
&&\int_0^1|f(\alpha)|^2d\alpha\\
&=&\sum\limits_{p_1}\sum\limits_{p_2}p_1^{1-\gamma}p_2^{1-\gamma}\int_0^1e(\alpha(p_1-p_2))d\alpha\\
&=&\mathop{\mathop{\sum\limits_{p\in\mathcal{A}_c}}}_{|p-\frac{n}{3}|<y}p^{2-2\gamma}\ll n^{2-2\gamma}yn^{\gamma-1}=n^{1-\gamma}y.
\end{eqnarray*}
\end{lemma}
\section{Proof of the theorem}
\begin{align*}
R(n)=&\mathop{\mathop{\mathop{\sum\limits_{p_1+p_2+p_3=n}}_{p_i\in \mathcal{A}_c}}_{|p_i-\frac{n}{3}|<y}}1=\int_0^1f^3(\alpha)e(-n\alpha)d\alpha\\
=&\int_0^1g^3(\alpha)e(-n\alpha)d\alpha+O\left(\int_0^1f^2(\alpha)(f(\alpha)-g(\alpha))d\alpha\right)\\&+O\left(\int_0^1f(\alpha)g(\alpha)(f(\alpha)-g(\alpha))d\alpha\right)+O\left(\int_0^1g^2(\alpha)(f(\alpha)-g(\alpha))d\alpha\right).
\end{align*}
By Cauchy inequality and trivial estimate, we have
\begin{align*}
R(n)=&\int_0^1g^3(\alpha)e(-n\alpha)d\alpha+O\left(\max\limits_{\alpha\in[0,1)}|f(\alpha)-g(\alpha)|\int_0^1|f^2(\alpha)|d\alpha\right)\\
&+O\left(\max\limits_{\alpha\in[0,1)}|f(\alpha)-g(\alpha)|\int_0^1|g^2(\alpha)|d\alpha\right).
\end{align*}
By Lemma 3.3, Lemma 3.4 and Lemma 3.5, we have
\begin{eqnarray*}
R(n)=\mathfrak{S}_3(n)\frac{3y^2}{\log^3 n}+O(\frac{y^2}{\log^4 n}),
\end{eqnarray*}
where $\mathfrak{S}_3(n)=\prod\limits_{p|n}(1-\frac{1}{(p-1)^2})\prod\limits_{p\nmid n}(1+\frac{1}{(p-1)^3})$. This complete the proof of the theorem.
\section*{Acknowledgements}
The author would thank his supervisor for his help and encouragement.


\end{document}